\documentclass[letterpaper,10pt,conference]{ieeeconf}
\IEEEoverridecommandlockouts
\usepackage[T1]{fontenc}
\usepackage[utf8]{inputenc}
\usepackage{amsmath,amssymb}
\usepackage{optidef}
\usepackage{algorithm}
\usepackage{color}
\usepackage{comment}
\usepackage{algpseudocode}
\usepackage{cite,url}
\usepackage{tikz}
\usetikzlibrary{arrows.meta}
\usepackage[hidelinks]{hyperref}

\usepackage{amsthm}
\theoremstyle{plain}
\newtheorem{assumption}{Assumption}
\newtheorem{theorem}{Theorem}
\newtheorem{proposition}{Proposition}
\newtheorem{lemma}{Lemma}

\newcommand{\R}{\mathbb R}
\newcommand{\Y}{\mathcal Y}
\newcommand{\cP}{\mathcal P}
\newcommand{\norm}[1]{\lVert#1\rVert}
\newcommand{\NLP}{\mathrm{NLP}}
\newcommand{\QP}{\mathrm{QP}}
\newcommand{\LP}{\mathrm{LP}}
\newcommand{\CIA}{\mathrm{CIA}}
\DeclareMathOperator*{\argmin}{arg\,min}

\title{\LARGE\bf Real-Time Nonlinear Mixed-Integer Model Predictive Control With Horizon-Tail Relaxations}

\author{Armin Nurkanovi\'c$^{1}$%
\thanks{*This research was supported by DFG via projects 504452366 (SPP 2364), 560056112 (robust MPC), 535860958 (ALeSCo) and 525018088 (MAWERO), and by BMWK via 03EN3054B.}%
\thanks{$^{1}$Systems Control and Optimization Laboratory, Department of Microsystems Engineering (IMTEK), University of Freiburg, Germany, {\tt\small armin.nurkanovic@imtek.uni-freiburg.de}}
\thanks{The author thanks Nathalie Janes and Moritz Diehl for interesting discussions on mixed-integer optimal control.
	OpenAI's Codex was used to assist with drafting and language editing, as well as implementing the accompanying code.
	The author reviewed and verified all AI-assisted content and takes full responsibility for the final manuscript and results.}
}

\begin{document}
\maketitle

\begin{abstract}
	Model predictive control (MPC) of nonlinear systems with discrete decisions leads to mixed-integer nonlinear programs (MINLPs) whose combinatorial complexity often prevents real-time deployment.
	In this paper, we propose a decomposition-based MPC method that reorders and reduces computations to shorten the delay between receiving a state measurement and applying the control input.
	Before the measurement arrives, we solve a continuous relaxation of the MINLP at a predicted state and construct a mixed-integer quadratic program (MIQP).
	Once the measured state is available, we solve the MIQP that reconstructs the discrete controls and simultaneously performs one sequential quadratic programming step toward the corresponding nonlinear solution.
	To reduce the combinatorial complexity in the horizon length, we either fix the integer controls or relax their integrality over the prediction horizon tail, while retaining an integer-feasible first control input.
	We derive error bounds relating the nonlinear objective values obtained with the horizon-tail approximations to the optimal value of the original mixed-integer problem.
	In a data center cooling example, the proposed method reduces mean feedback computation time by $94\%$ relative to a state-of-the-art MINLP method with a total closed-loop cost difference of less than $1\%$.
\end{abstract}

\section{Introduction}\label{sec:introduction}
Model predictive control (MPC) is a feedback control strategy that repeatedly solves optimal control problems.
Its appeal lies in its ability to directly specify objectives and constraints while using predictive models to anticipate the effects of control decisions.
MPC algorithms are well established for linear and nonlinear systems with continuous control inputs~\cite{Diehl2009c}.
Yet, many real-world systems also require discrete controls.

Combined with nonlinear dynamics, these discrete decisions give rise to nonlinear mixed-integer optimal control problems (MIOCPs), whose time discretization yields mixed-integer nonlinear programs (MINLPs).
Generic MINLP solvers~\cite{Belotti2013}, even when terminated early, are usually too slow for real-time MPC deployment.
Particularly for nonlinear systems, this has motivated methods tailored to MIOCPs that rapidly compute high-quality approximate solutions, often with well-understood theoretical properties.
These methods include, for example, switching-time optimization~\cite{Gerdts2006,Abbasi2026}, decomposition methods~\cite{Sager2011a,Kirches2020} or machine-learning-based approaches~\cite{Orrico2026,Boldocky2026,Adamek2026}.

Decomposition methods for MIOCPs solve continuous nonlinear programs (NLPs) related to the original MINLP, combined with a discrete reconstruction step~\cite{Sager2011a,Sager2009,Kirches2020}, yielding high-quality approximations with well-understood error bounds~\cite{Sager2011a,Kirches2020}.
Discrete reconstruction recovers the integer controls using different methods ranging from OCP-tailored rounding~\cite{Sager2009} and variants accommodating dwell-time constraints~\cite{Zeile2021,Chen2022} to integer programming~\cite{Sager2011a,Buerger2020a} and tailored shortest-path methods incorporating combinatorial constraints~\cite{Bestehorn2021}.
However, recovering integer controls without explicitly accounting for their effects on the dynamics can cause performance losses or infeasibility on coarse time-discretization grids~\cite{Buerger2020b}.
Reconstruction through a mixed-integer quadratic program (MIQP) incorporates local approximations of the dynamics and objective, improving solution quality at greater computational cost~\cite{Buerger2023}.
For real-time MPC, this trade-off motivates further approximations that reduce computation while preserving solution quality, a direction pursued in this paper.

Optimization problem solution approximations are already standard in real-time MPC with continuous controls, where their theoretical properties are well understood~\cite{Rawlings2017}.
Rawlings and Risbeck~\cite{Rawlings2017a} show how these properties extend conceptually to approximate mixed-integer MPC.
This motivates us to adapt algorithmic ideas from real-time continuous MPC to the mixed-integer setting~\cite{Diehl2009c}.

A shared idea in many of these algorithms is to reorder optimization computations to reduce the feedback delay~\cite{Diehl2009c}.
The computations are split into a preparation phase, which performs operations before the next state measurement $\bar{x}_0$ is known, such as evaluating functions and derivatives or presolving related problems, and a feedback phase, which computes a suitable control approximation with minimal effort once $\bar{x}_0$ is known.
Advanced-step methods compute in the preparation phase an (approximate) solution of the MPC problem at a predicted state $\tilde{x}_0$~\cite{Zavala2009,Nurkanovic2019a,Frey2024a}.

In this paper, we adapt these ideas to mixed-integer MPC by reordering and reducing the computations of decomposition algorithms with.
In the preparation phase, we solve the continuous relaxation of the MINLP, and optionally an MIQP at a predicted advanced state $\tilde{x}_0$.
In the feedback phase, we update the prepared MIQP with the measured state $\bar{x}_0$ and solve it using the available warm start.
To reduce the MIQP's combinatorial complexity, we either fix the integer controls or relax their integrality over the prediction horizon tail.
This reduces the feedback delay while providing bounds on the resulting performance losses.


\textit{Contributions.} The main contributions of this paper are: 
(i) a real-time mixed-integer nonlinear MPC method that separates preparation and feedback computations, combining integer reconstruction and one SQP step in a single MIQP; 
(ii) two horizon-tail approximations that fix integer controls or relax their integrality beyond a selected horizon, reducing the combinatorial search while retaining an integer-feasible first control input; 
(iii) bounds relating the nonlinear objective values obtained with these approximations to the globally optimal value of the original MINLP; 
(iv) a numerical comparison with a state-of-the-art MINLP method and a rule-based controller on a data center cooling example; and 
(v) an open-source implementation supporting several algorithmic variants and several MPC examples from the literature.\footnote{\url{https://github.com/nosnoc/advanced-minmpc}}

\section{Decomposition methods for MIOCP}\label{sec:decomposition}
In this section, we state the MIOCP of interest, and recall a standard decomposition method for solving such problems.
\subsection{Problem formulation}
We consider a discrete-time nonlinear system with states $x_k\in\R^{n_x}$, continuous controls $u_k\in\R^{n_u}$, and binary controls $v_k\in\{0,1\}^{n_v}$.
General bounded integer controls can be encoded by binary variables, so this restriction is without loss of generality.
For an initial state $\bar x_0$ and a horizon of $N$ intervals, define $x=(x_0,\ldots,x_N)$, $u=(u_0,\ldots,u_{N-1})$, and $v=(v_0,\ldots,v_{N-1})$, to be the discrete-time states and controls.
The MIOCP is stated as follows:
\begin{mini!}|s|[2]
 {x,u,v}{\displaystyle\sum_{k=0}^{N-1}\ell_k(x_k,u_k,v_k)+\ell_N(x_N)\label{eq:ocp-objective}}
 {\label{eq:ocp}}{}
 \addConstraint{x_0}{= \bar x_0,\label{eq:ocp-initial}}
 \addConstraint{x_{k+1}}{=f_k(x_k,u_k,v_k),\label{eq:ocp-dynamics}}{\quad k=0,\ldots,N-1}
 \addConstraint{h_k(x_k,u_k,v_k)}{\leq0\label{eq:ocp-path},}{\quad k=0,\ldots,N-1}
 \addConstraint{v_k}{\in\{0,1\}^{n_v}\label{eq:ocp-binary},}{\quad k=0,\ldots,N-1}
 \addConstraint{v}{\in\cP.\label{eq:ocp-combinatorial}}
\end{mini!}
Here, $\ell_k:\mathbb{R}^{n_x}\times\mathbb{R}^{n_u}\times\{0,1\}^{n_v}\to\R$ is the stage cost and $\ell_N:\R^{n_x}\to\R$ is the terminal cost.
The map $f_k:\R^{n_x}\times\R^{n_u}\times\{0,1\}^{n_v}\to\R^{n_x}$ denotes the discrete-time state-transition map.
The functions $h_k:\R^{n_x}\times\R^{n_u}\times \{0,1\}^{n_v}\to\R^{n_{h_k}}$ collect the path inequality constraints.
The polyhedron $\cP\subseteq[0,1]^{Nn_v}$ contains additional combinatorial constraints coupling several control intervals, e.g., dwell-time constraints~\cite{Chen2022}.

Next, we define a few shorthands.
We aggregate the binary controls as $y:=v\in\{0,1\}^{n_y}$, where $n_y=Nn_v$, and collect the continuous primal variables into $w:=(x,u)\in\R^{n_w}$.
Including the bounds in its linear representation, we write $\cP=\{y\in\R^{n_y}\mid Ay\leq b\}$.
The feasible binary set is $ \Y:=\cP\cap\{0,1\}^{n_y}$. 
The vector $p \in \R^{n_p}$ denotes the parameter, in our case, without loss of generality, this is the initial state $p = \bar x_0$.
Let $G$ collect the equality constraints in~\eqref{eq:ocp-initial} and~\eqref{eq:ocp-dynamics}, and let $H$ collect the inequality constraints in~\eqref{eq:ocp-path}.
A function $F$ is a generic nonlinear function collecting the stage and terminal costs in~\eqref{eq:ocp}.
In the remainder of this and the next section, we assume that the parameter $p$ is fixed.
Using this notation, problem~\eqref{eq:ocp} can be written compactly as
\begin{mini!}|s|[2]
 {y,w}{F(y,w)\label{eq:minlp-objective}}
 {\label{eq:minlp}}{}
 \addConstraint{G(y,w;p)}{=0\label{eq:minlp-equality}}
 \addConstraint{H(y,w)}{\leq0\label{eq:minlp-inequality}}
 \addConstraint{y}{\in\Y.\label{eq:minlp-integer}}
\end{mini!}

It will be later useful to define a standard NLP where $y \in \mathcal P$ is a parameter and not a decision variable:
\begin{mini!}|s|[2]
	{w}{F(y,w)\label{eq:nlp-objective}}
	{\label{eq:nlp}}{J_{\NLP}(y):=}
	\addConstraint{G(y,w;p)}{=0\label{eq:nlp-equality}}
	\addConstraint{H(y,w)}{\leq0 .\label{eq:nlp-inequality}}
\end{mini!}
Here, $J_{\NLP}(y)$ denotes the optimal value function of this NLP.
Note that we can compactly express the integer part of the MINLP~\eqref{eq:minlp} solution as~\cite{Buerger2023}:
\begin{equation}\label{eq:compact-minlp}
	y^\circ =  \mathrm{arg}\min_ y \ J_{\NLP}(y) \ \mathrm{s.t.} \  y \in \Y.
\end{equation}
A relaxed integer solution $y^*$, and thus a lower bound for $J_{\NLP}(y^*) \leq J_{\NLP}(y^\circ)$ is obtained by solving the relaxed MINLP~\eqref{eq:minlp}, compactly written as:
$ y^* =  \mathrm{arg}\min_ y \ J_{\NLP}(y) \ \mathrm{s.t.} \  y \in \cP$.

\subsection{A three-phase decomposition method}\label{sec:three_phase}
For optimal control and MPC applications, decomposition methods typically compute high-quality approximate solutions of~\eqref{eq:minlp} more efficiently than generic MINLP methods~\cite{Sager2011a,Sager2009,Kirches2020}.
The decomposition method introduced by Sager~\cite{Sager2009} consists of three steps: solve a relaxation of the MINLP, integer reconstruction, and nonlinear polishing.
Following the exposition in~\cite{Buerger2023}, these steps are:
\begin{enumerate}
		\renewcommand{\theenumi}{S\arabic{enumi}}
	\renewcommand{\labelenumi}{(\theenumi)}
	\item\label{step:relaxation} The integrality condition in~\eqref{eq:minlp} is relaxed by replacing $\Y$ with $\cP\subseteq\R^{n_y}$, yielding a nonlinear program (NLP).
	A local minimizer of this NLP is denoted by $(y^*,w^*)$, where $y^*\in\cP$ but not necessarily $y^*\in\Y$.
	If $y^* \in \Y$, then steps two and three are omitted.
	\item\label{step:reconstruction} An integer solution $\hat y\in\Y$ is reconstructed by minimizing a distance function:
	\begin{align}\label{eq:generic-distance}
		\hat y \in \argmin_{y\in\Y} d(y;y^*).
	\end{align}
	\item\label{step:polishing} The integer variables are fixed to $y=\hat y$, and the resulting NLP~\eqref{eq:nlp}, referred to as the polishing NLP, is solved.
	Whenever this NLP admits a solution $w^*(\hat y)$, the pair $(\hat y,w^*(\hat y))$ is feasible for~\eqref{eq:minlp}.
\end{enumerate}

This procedure typically yields a pair $(\hat y,w^*(\hat y))$ that is locally optimal for the fixed-integer NLP~\eqref{eq:nlp} and feasible, but generally suboptimal, for the MINLP~\eqref{eq:minlp}.
Let $(y^\circ,w^\circ)$ denote a global minimizer of~\eqref{eq:minlp}.
If the continuous relaxation is solved globally, the optimal values satisfy $J_{\NLP}(y^*)\leq J_{\NLP}(y^\circ)\leq J_{\NLP}(\hat y)$~\cite{Buerger2023}.
Any local minimizer $(y^*,w^*)$ of the relaxed NLP can also be used for integer reconstruction in \eqref{step:reconstruction}, although the first inequality does not necessarily hold anymore.
The practical performance and computational complexity of decomposition methods depend strongly on the choice of the distance function $d(y;y^*)$.
In the remainder of the section, we review existing choices for $d$.
In Sec.~\ref{sec:relaxations} we propose new ones tailored to MPC, and establish their relations to $J_{\NLP}(y^\circ)$.

\subsubsection{Ideal distance function}
The ideal choice for the distance function is $d_{\NLP}(y;y^*):=J_{\NLP}(y)-J_{\NLP}(y^*)$.
If $y^*$ is a global solution of the continuous relaxation, this function is nonnegative for all $y\in\Y$, since $J_\mathrm{NLP}(y^*)\leq J_\mathrm{NLP}(y)$.
Because $J_\mathrm{NLP}(y^*)$ is constant, solving~\eqref{eq:generic-distance} with this choice is equivalent to solving~\eqref{eq:compact-minlp}, and hence the original MINLP~\eqref{eq:minlp}.
This choice is therefore usually intractable and we look for more tractable approximations of $d_{\NLP}(y;y^*)$.


\subsubsection{Combinatorial integral approximation (CIA) distance}
A simpler choice is $d(y;y^*):=\|y-y^*\|_q$, which yields a mixed-integer linear programming (MILP) formulation of~\eqref{eq:generic-distance} for $q=1$ or $q=\infty$.
The frequently used CIA norm $d_{\CIA}(y;y^*):=\|U(y-y^*)\|_\infty$, with an invertible matrix $U$ defined in~\cite{Sager2011a}, also admits an MILP formulation.
If $\mathcal P = \R^{n_y}$, sum-up rounding provides a polynomial-time solution to the CIA problem under suitable assumptions~\cite[Theorem~3]{Sager2011a}.
Norm-based distance functions are attractive for their computational efficiency, as they admit tailored solution methods, such as those implemented in \texttt{pycombina}~\cite{Buerger2020a}, which we use in this work.
On the downside, such approximations reconstruct $\hat y$ based solely on $y^*$, disregarding $w^*$ or the nonlinear dynamics and constraints, which can lead to performance losses~\cite{Buerger2020b}.

\subsubsection{Quadratic-model distance}
Since minimizing $d(y;y^*):=J_{\NLP}(y)-J_{\NLP}(y^*)$ over $\Y$ is computationally intractable, \cite{Buerger2023} proposes a quadratic approximation of $J_{\NLP}(y)$, turning~\eqref{eq:generic-distance} into an MIQP.
We recall this approximation here and further adapt it in Sec.~\ref{sec:relaxations}.

Define the Lagrangian of the fixed-$y$ NLP~\eqref{eq:nlp} as
\begin{equation*}
	\mathcal L(y,w,\lambda,\mu;p)
	:=F(y,w)+\lambda^\top G(y,w;p)+\mu^\top H(y,w).
\end{equation*}
The Lagrangian of the continuous relaxation of~\eqref{eq:minlp} is
$\mathcal L_{\mathrm{rel}}(y,w,\lambda,\mu,\rho;p)
:=\mathcal L(y,w,\lambda,\mu;p)+\rho^\top(Ay-b)$.
Since the additional term is affine in $y$, both Lagrangians have the same Hessian with respect to $(y,w)$.
At a Karush--Kuhn--Tucker (KKT) point $(y^*,w^*,\lambda^*,\mu^*,\rho^*)$ of the relaxed NLP in \eqref{step:relaxation}, we approximate this Hessian by a symmetric positive definite matrix
$B:=B(y^*,w^*,\lambda^*,\mu^*)\succ0$.
With $\Delta y:=y-y^*$ and $\Delta w:=w-w^*$, define the linearization of $F$ compactly as
\begin{align*}
	F_{\mathrm L}(y,w;y^*,w^*):=F(y^*,w^*)+
	\nabla_{(y,w)}F(y^*,w^*)^\top
	\begin{bmatrix}
		\Delta y\\ \Delta w
	\end{bmatrix}\!.
\end{align*}
The linearizations $G_{\mathrm L}(y,w;p,y^*,w^*)$ and $H_{\mathrm L}(y,w;y^*,w^*)$ are defined analogously.
Furthermore, a quadratic approximation of the objective of~\eqref{eq:minlp} using $B \succ 0$ is defined as
\begin{align*}
	F_{\QP}(y,w;y^*\!,\!w^*,\!B)\! :=\! F_{\mathrm L}(y,w;y^*,\!w^*)\!+\!\frac12
	\begin{bmatrix}\Delta y\\\Delta w\end{bmatrix}^{\!\top}
	\!\!\! \! B\begin{bmatrix}\Delta y\\\Delta w\end{bmatrix}\!.
\end{align*}
For the linearization point $(y^*,w^*)$, the MIQP approximation of~\eqref{eq:minlp} reads
\begin{mini!}|s|[2]
	{y,w}{F_{\QP}(y,w;y^*,w^*,B)\label{eq:miqp-objective}}
	{\label{eq:miqp}}{}
	\addConstraint{G_{\mathrm L}(y,w;p,y^*,w^*)}{=0\label{eq:miqp-equality}}
	\addConstraint{H_{\mathrm L}(y,w;y^*,w^*)}{\leq0\label{eq:miqp-inequality}}
	\addConstraint{y}{\in\Y.\label{eq:miqp-combinatorial}}
\end{mini!}
In analogy to~\eqref{eq:nlp}, for fixed $y$, the quadratic value function is defined as
\begin{mini!}|s|[2]
	{w}{F_{\QP}(y,w;y^*,w^*,B)\label{eq:qp-value-objective}}
	{\label{eq:qp-value}}{J_{\QP}(y;y^*,w^*,B):=}
	\addConstraint{G_{\mathrm L}(y,w;p,y^*,w^*)}{=0\label{eq:qp-value-equality}}
	\addConstraint{H_{\mathrm L}(y,w;y^*,w^*)}{\leq0.\label{eq:qp-value-inequality}}
\end{mini!}
For $B\succ0$, $J_{\QP}$ is convex, continuous, and piecewise quadratic in $y$ on its feasible domain~\cite{Buerger2023}.
Moreover, $J_{\QP}(y;y^*,w^*,B)$ approximates $J_{\NLP}(y)$, with an error bound given in Lemma~\ref{lem:value-consistency} below.
For fixed $p$, the quadratic-model distance is defined as
\begin{equation}\label{eq:qp-distance}
	d_{\QP}(y;y^*,w^*,B):=J_{\QP}(y;y^*,w^*,B)-J_{\NLP}(y^*).
\end{equation}
The relaxed NLP's KKT point $(y^*,w^*)$ also satisfies the KKT conditions of the MIQP's~\eqref{eq:miqp} continuous relaxation, which is convex because $B\succ0$ and its constraints are affine.
Hence, $(y^*,w^*)$ globally minimizes this QP, and fixing any feasible $y\in\cP$ cannot decrease its optimal value, so $J_{\QP}(y)\geq F(y^*,w^*)\geq J_{\NLP}(y^*)$, hence $d_{\QP}(y;y^*,w^*,B) \geq0$.
Minimizing this distance over $\Y$ is equivalent to solving~\eqref{eq:miqp} constructed at $(y^*,w^*,p)$.
Note that we could use also $B = 0$, then~\eqref{eq:miqp} reduces to an MILP.

\section{Relaxations in the integer reconstruction}\label{sec:relaxations}
Solving the MIQP~\eqref{eq:miqp} is usually the most expensive step of the decomposition approach with a quadratic-model distance.
Since MPC applies only the first control input, only $v_0$ must be integer feasible.
This motivates reducing the combinatorial complexity through two horizon-tail approximations inspired by~\cite{Frick2015a}.
While~\cite{Frick2015a} uses these approximations in a tailored solver for fixed-matrix MIQPs, we study them within the decomposition approach and apply them to real-time MPC in Sec.~\ref{sec:advanced}.

For either approximation, partition $y=(y_1,y_2)$ at a horizon $0\leq M\leq N$, where
\begin{equation}\label{eq:horizon-blocks}
	y_1=(v_0,\ldots,v_{M-1}),\qquad y_2=(v_M,\ldots,v_{N-1}).
\end{equation}
The corresponding block is empty when $M=0$ or $M=N$.
The objective and linearized constraints retain the full horizon.

\begin{figure}[t]
	\centering
	\includegraphics[width=\columnwidth]{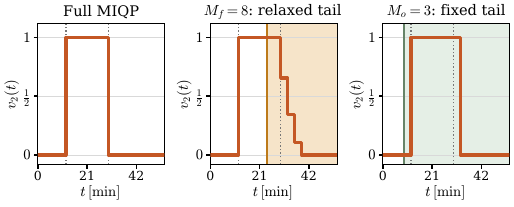}
	\vspace{-0.8cm}
	\caption{Horizon-tail approximations for a binary control in the parallel-chiller example of Sec.~\ref{sec:example}.
		The full MIQP optimizes all binary decisions.
		The feasibility relaxation ($M_f=8$) allows a continuous tail, whereas the optimality relaxation ($M_o=3$) fixes the tail to prescribed binary values.
		Shading indicates the respective horizon tails.}
	\label{fig:relaxed-horizons}
	\vspace{-0.6cm}
\end{figure}

\subsection{Relaxation of optimality}
Suppose we have an integer-feasible vector $y^{**}\in\Y$.
In practice, it may be obtained via sum-up rounding or by solving the CIA-MILP with \texttt{pycombina}~\cite{Buerger2020a}.
For $0<M_o\leq N$, the first block remains a binary decision and the second block is fixed to the corresponding part of $y^{**}$.
Thus, in \eqref{step:reconstruction} we solve an MIQP with fewer free binary variables:
\vspace{-0.2cm}
\begin{mini!}|s|[2]
	{y_1,y_2,w}{F_{\QP}(y,w;y^*,w^*,B)\label{eq:optimality-objective}}
	{\label{eq:optimality}}{}
	\addConstraint{G_{\mathrm L}(y,w;p,y^*,w^*)}{=0\label{eq:optimality-equality}}
	\addConstraint{H_{\mathrm L}(y,w;y^*,w^*)}{\leq0\label{eq:optimality-inequality}}
	\addConstraint{y}{\in\cP\label{eq:optimality-combinatorial}}
	\addConstraint{y_1}{\in\{0,1\}^{M_on_v}\label{eq:optimality-binary}}
	\addConstraint{y_2}{=y_2^{**}.\label{eq:optimality-tail}}
\end{mini!}
Every entry of $y$ is binary, so~\eqref{eq:optimality} supplies an integer-feasible first input.
For $M_o=0$, the complete integer sequence is fixed to $y^{**}$ and only one continuous QP is solved.
For $M_o=N$, problem~\eqref{eq:optimality} is the full MIQP~\eqref{eq:miqp}.
For any $M_o$, subsequent nonlinear polishing yields a feasible, potentially suboptimal solution of the full MINLP~\eqref{eq:minlp}, provided the polishing NLP in \eqref{step:polishing} is feasible.
\subsection{Relaxation of feasibility}
For $1\leq M_f\leq N$, only the first block is required to be binary, while integrality is relaxed in the remaining block.
We solve the following MIQP:
\allowdisplaybreaks
\begin{mini!}|s|[2]
	{y_1,y_2,w}{F_{\QP}(y,w;y^*,w^*,B)\label{eq:feasibility-objective}}
	{\label{eq:feasibility}}{}
	\addConstraint{G_{\mathrm L}(y,w;p,y^*,w^*)}{=0\label{eq:feasibility-equality}}
	\addConstraint{H_{\mathrm L}(y,w;y^*,w^*)}{\leq0\label{eq:feasibility-inequality}}
	\addConstraint{y}{\in\cP\label{eq:feasibility-combinatorial}}
	\addConstraint{y_1}{\in\{0,1\}^{M_fn_v}\label{eq:feasibility-binary}}
	\addConstraint{y_2}{\in[0,1]^{(N-M_f)n_v}.\label{eq:feasibility-tail}}
\end{mini!}
Since $M_f\geq1$, the first control input is integer feasible, while the second block can be fractional.
For $M_f=N$, this problem is the full MIQP~\eqref{eq:miqp}.
With this relaxation, the decomposition no longer necessarily provides a feasible solution of the full MINLP~\eqref{eq:minlp}.
Nevertheless, the relaxed problem provides a lower bound on the full MIQP value, as discussed below, and retains an integer-feasible first input for MPC.
From a dynamic programming perspective, both approximations replace the full MIQP's tail value function by one obtained with either prescribed binary controls or relaxed integrality.
Fig.~\ref{fig:relaxed-horizons} illustrates both approximations and compares them with the full MIQP.

\subsection{Relation to the nonlinear value function}
Throughout this section, $p$ is fixed and we abbreviate $J_{\QP}(y):=J_{\QP}(y;y^*,w^*,B)$.
We now relate the optimal value functions of the relaxed MIQPs~\eqref{eq:optimality} and~\eqref{eq:feasibility} to $J_{\NLP}(y)$, the optimal value function of~\eqref{eq:nlp}.
We begin by relating the optimal value functions of the three MIQP variants: the full MIQP and its optimality and feasibility relaxations.
\begin{theorem}[Horizon-tail gap bound]\label{th:horizon-gap}
	Suppose that the continuous relaxation of~\eqref{eq:miqp} has a compact feasible set and fix $(p,y^*,w^*,B)$ with $B\succ0$.
	Fix a binary sequence $y^{**}$ feasible for~\eqref{eq:miqp} and assume that, for every partition~\eqref{eq:horizon-blocks}, replacing the tail of any feasible sequence with a binary prefix by $y_2^{**}$ preserves feasibility after reoptimizing $w$.
	Let $J_{\QP}^*$, $U_M$, and $L_M$ denote the globally optimal values of~\eqref{eq:miqp}, \eqref{eq:optimality}, and~\eqref{eq:feasibility}, respectively, with $M_o=M_f=M$ and $1\leq M\leq N$.
	Then there exists $K>0$, independent of $M$, such that
	\begin{subequations}\label{eq:horizon-gap}
		\begin{align}
			&L_M\leq J_{\QP}^*\leq U_M, \label{eq:horizon-gap-opt}\\
			&0\leq U_M-L_M\leq K\sqrt{(N-M)n_v}. \label{eq:horizon-gap-M}
		\end{align}
	\end{subequations}
	Moreover, $L_M$ is nondecreasing, $U_M$ is nonincreasing, and $L_N=U_N=J_{\QP}^*$.
\end{theorem}

\begin{proof}
	Every feasible point of the full MIQP~\eqref{eq:miqp} is feasible for the feasibility-relaxed problem~\eqref{eq:feasibility}, hence $L_M\leq J_{\QP}^*$.
	Conversely, every optimizer of~\eqref{eq:optimality} is feasible for the full MIQP~\eqref{eq:miqp}, hence $J_{\QP}^*\leq U_M$.
	Together, these inequalities establish~\eqref{eq:horizon-gap-opt}.
	Increasing $M$ shrinks the feasible set of~\eqref{eq:feasibility} and enlarges that of~\eqref{eq:optimality}, proving monotonicity of $L_M$ and $U_M$.
	At $M=N$, both problems coincide with~\eqref{eq:miqp}.
	
	By~\cite[Theorem~6.7]{Borrelli2017}, $J_{\QP}(y)$ is continuous and piecewise quadratic on its feasible domain.
	Its finitely many quadratic pieces and compactness of the feasible set of the relaxation of~\eqref{eq:miqp} imply Lipschitz continuity of $J_{\QP}(y)$ with a constant $K>0$.
	Let $\hat y_f=(\hat y_{f,1},\hat y_{f,2})$ be the $y$-component of a global optimizer of~\eqref{eq:feasibility}, so that $L_M=J_{\QP}(\hat y_f)$.
	Similarly, let $(\hat y_{o,1},y_2^{**})$ be the corresponding component of a global optimizer of~\eqref{eq:optimality}, so that $U_M=J_{\QP}(\hat y_{o,1},y_2^{**})$.
	Then
	\begin{align*}
		U_M-L_M
		&=J_{\QP}(\hat y_{o,1},y_2^{**})-J_{\QP}(\hat y_{f,1},\hat y_{f,2})\\
		&\leq J_{\QP}(\hat y_{f,1},y_2^{**})-J_{\QP}(\hat y_{f,1},\hat y_{f,2})\\
		&\leq K\norm{y_2^{**}-\hat y_{f,2}}
		\leq K\sqrt{(N-M)n_v}.
	\end{align*}
	In the first inequality, we used $J_{\QP}(\hat y_{o,1},y_2^{**})\leq J_{\QP}(\hat y_{f,1},y_2^{**})$, which follows because the tail-replacement assumption makes $(\hat y_{f,1},y_2^{**})$ feasible in~\eqref{eq:optimality}, while $\hat y_{o,1}$ is optimal.
	The second inequality follows from Lipschitz continuity of $J_{\QP}$. 
	The last uses the fact that $y_2$ has $q=n_v(N-M)$ components, each differing by at most one.
	Thus, $\norm{y_2^{**}-\hat y_{f,2}}\leq\sqrt{\sum_{i=1}^{q}1^2}=\sqrt{q}$, which establishes~\eqref{eq:horizon-gap-M}.
\end{proof}
We proceed towards relating $J_{\NLP}(y)$ to $J_{\QP}(y)$.
Recall that $y$ is interpreted here as a fixed parameter in the NLP~\eqref{eq:nlp}, for which we will assume some standard regularity conditions.
\begin{assumption}[Local NLP regularity]\label{ass:nlp_regularity}
	Let $(w^*,\lambda^*,\mu^*)$ be a KKT point of~\eqref{eq:nlp} at $y=y^*$, and define $\mathcal A^*:=\{i\mid H_i(y^*,w^*)=0\}$ and $\mathcal A_+^*:=\{i\in\mathcal A^*\mid\mu_i^*>0\}$.
	It holds that:
	\begin{enumerate}
		\item[(a)] The functions $F$, $G$, and $H$ are twice continuously differentiable near $(y^*,w^*)$.
		\item[(b)] Linear independence constraint qualification (LICQ) holds with respect to $w$, i.e., the vectors
		\begin{align*}
		\{\nabla_wG_j(y^*,w^*;p)\}_j
		\cup\{\nabla_wH_i(y^*,w^*)\}_{i\in\mathcal A^*}	
		\end{align*}
		are linearly independent.
		\item[(c)] Strong second order sufficient condition (SSOSC) holds with respect to $w$, i.e.,
		\begin{align*}
		d_w^\top\nabla_{ww}^2\mathcal L(y^*,w^*,\lambda^*,\mu^*;p)d_w>0	
		\end{align*}
		for every nonzero $d_w$ in the strong critical cone
		\begin{equation*}
					\mathcal C^{\mathrm S}:=
			\left\{d_w\ \middle|\
			\begin{aligned}
				\nabla_wG(y^*,w^*;p)^\top d_w&=0,\\
				\nabla_wH_i(y^*,w^*)^\top d_w&=0,\quad i\in\mathcal A_+^*
			\end{aligned}\right\}.
		\end{equation*}
	\end{enumerate}
\end{assumption}


\vspace{0.2cm}
Next, we establish an intermediate result similar to~\cite[Corollary~1]{Buerger2023}, but more directly using only elementary NLP perturbation theory.
Throughout the following analysis, $w_{\NLP}(y)$ denotes the locally unique NLP solution map in the polishing step used in Lemma~\ref{lem:value-consistency}, with $w_{\NLP}(y^*)=w^*$ and $J_{\NLP}(y):=F(y,w_{\NLP}(y))$.
\begin{lemma}\label{lem:value-consistency}
	Suppose Assumption~\ref{ass:nlp_regularity} holds for the NLP~\eqref{eq:nlp} at the fixed parameter $y^*$ and the corresponding solution $w^*$, and that $B\succ0$ in~\eqref{eq:qp-value}.
		Then there exist $\epsilon>0$ and $C\geq0$ such that for all $y \in \mathcal{B}_\epsilon(y^*):= \{ y \in \R^{n_y} \mid \norm{y-y^*}<\epsilon\}$:
	\begin{equation}\label{eq:value-consistency}
			|J_{\NLP}(y)-J_{\QP}(y)|\leq C\norm{y-y^*}^2.
		\end{equation}
\end{lemma}
\begin{proof}
	Fix the parameter $y=y^*$ also in the QP~\eqref{eq:qp-value}.
	By construction, the NLP and QP have the same KKT conditions at $(w^*,\lambda^*,\mu^*)$.
	Since $B\succ0$, the unique solution of this QP is therefore $w^*$, with $\Delta w=0$.
	Hence, $J_{\QP}(y^*)=J_{\NLP}(y^*)$, and it can be verified by direct computation that $\nabla_y\mathcal L(y^*,w^*,\lambda^*,\mu^*;p)	=\nabla_y\mathcal L_{\QP}(y^*,w^*,\lambda^*,\mu^*;p)$, where $\mathcal L_{\QP}$ is the Lagrangian of the QP~\eqref{eq:qp-value}.
	
	By Assumption~\ref{ass:nlp_regularity} and~\cite[Theorem~2.4.5~(iv)]{Fiacco1983}, the local primal--dual solution of the NLP~\eqref{eq:nlp} is Lipschitz continuous in $y$.
	Because $B\succ0$ ensures SSOSC and Assumption~\ref{ass:nlp_regularity} implies LICQ also for the QP~\eqref{eq:qp-value} at $w^*$, we can apply~\cite[Theorem~2.4.5~(iv)]{Fiacco1983} to the QP's solution map as well.
	Thus, let
	$
	z_i(y)=(w_i(y),\lambda_i(y),\mu_i(y)),
	\ i\in\{\NLP,\QP\},
	$
	denote these primal--dual solution maps, with $z_i(y^*)=(w^*,\lambda^*,\mu^*)$.
	By twice continuous differentiability of $F$, $G$, and $H$ (Ass.~\ref{ass:nlp_regularity}(a)) and Lipschitz continuity of the solution maps, $\nabla_y\mathcal L$ and $\nabla_y\mathcal L_{\QP}$, evaluated along the respective solution maps, are locally Lipschitz in $y$.

	Furthermore, by applying \cite[Theorem 2.4.5~(v)]{Fiacco1983} to the NLP and QP, we have that
	$\nabla_y J_{\NLP}(y^*)  = \nabla_y \mathcal L(y^*,w^*,\lambda^*,\mu^*;p)$
	and $\nabla_y J_{\QP}(y^*)  = \nabla_y \mathcal L_{\QP}(y^*,w^*,\lambda^*,\mu^*;p)$.
	The same identities hold locally along the respective solution maps, so both value-function gradients are locally Lipschitz.
	
	Finally, define $E(y):=J_{\NLP}(y)-J_{\QP}(y)$.
	It follows from the arguments above that there exists $\epsilon>0$ such that $\nabla E(y)$ is Lipschitz with a constant $L_E$ for $y \in \mathcal{B}_\epsilon(y^*)$.
	Let $d_y:=y-y^*$.
	Using $E(y^*)=0$ and $\nabla E(y^*)=0$ and applying the fundamental theorem of calculus for $\norm{d_y}<\epsilon$, with 
	$|E(y^* + d_y) + {E(y^*)}| = |E(y)|$ we obtain
		\begin{align*}
		&|E(y)|=\left| E(y^*) + \int_0^1[\nabla E(y^*+t {d_y})-\nabla E(y^*)]^\top {d_y}\, \mathrm{d} t\right|\\
		&\leq\int_{0}^{1} L_E t \norm{d_y}^2\,\mathrm{d}t
		=\frac{1}{2} L_E \norm{d_y}^2.
	\end{align*}
	Thus,~\eqref{eq:value-consistency} holds with $C=L_E/2$.
\end{proof}

Finally, we relate the value $J_{\NLP}(\hat y)$ obtained by the decomposition approach to the globally optimal MINLP value $J_{\NLP}(y^\circ)$ in~\eqref{eq:compact-minlp}.
Lemma~\ref{lem:value-consistency} requires $y\in\mathcal B_\epsilon(y^*)$, which need not hold for arbitrary binary vectors $y$ and a relaxed solution $y^*$.
Hence, we make an assumption specific to mixed-integer OCPs, analogous to~\cite[Assumption~1]{Buerger2023}, that allows us to apply the preceding result.
\begin{assumption}\label{ass:integer-locality}
	Let $\epsilon>0$ be as in Lemma~\ref{lem:value-consistency}, and fix $M\in\{1,\ldots,N\}$.
	Let $\hat y$ be the $y$-component of a global optimizer of either~\eqref{eq:optimality} with $M_o=M$ or~\eqref{eq:feasibility} with $M_f=M$.
	Assume that $\hat y\in\mathcal B_\epsilon(y^*)$ and that there exists a global minimizer $(y^\circ,w^\circ)$ of~\eqref{eq:minlp} satisfying
	$y^\circ\in\mathcal B_\epsilon(y^*)$, $w^\circ=w_{\NLP}(y^\circ)$.
\end{assumption}

Recall that $(y^*,w^*)$ is the relaxed solution obtained in \eqref{step:relaxation}, while $J_{\NLP}(\hat y)$ is the value of the polished solution obtained in \eqref{step:polishing} by fixing $y=\hat y$ in~\eqref{eq:nlp}, with $\hat w=w_{\NLP}(\hat y)$.
Next, our main theoretical result follows.

\begin{theorem}[Nonlinear value bounds]\label{th:nonlinear-loss}
	Fix $p$ and let $(y^*,w^*)$ be a local solution of the continuous relaxation of~\eqref{eq:minlp}.
	Suppose Assumption~\ref{ass:nlp_regularity} holds for the fixed-$y^*$ NLP~\eqref{eq:nlp} at $w^*$.
	Construct the MIQP~\eqref{eq:qp-value} at $(y^*,w^*)$ with $B\succ0$, and let $\epsilon$ and $C$ be given by Lemma~\ref{lem:value-consistency}.
	Suppose the Assumption of Theorem~\ref{th:horizon-gap} hold for this MIQP, and let $K$ be its gap-bound constant.
	
	Fix $M\in\{1,\ldots,N\}$ and let $\hat y$ be the $y$-component of a global optimizer of either~\eqref{eq:optimality} with $M_o=M$ or~\eqref{eq:feasibility} with $M_f=M$.
	Let $(y^\circ,w^\circ)$ be a global minimizer of~\eqref{eq:minlp}, and suppose Assumption~\ref{ass:integer-locality} holds for $\hat y$ and $(y^\circ,w^\circ)$ with the radius $\epsilon$.
	For the polishing solution $\hat w=w_{\NLP}(\hat y)$, define $J_{\NLP}(\hat y)=F(\hat y,\hat w)$; by the same assumption, $J_{\NLP}(y^\circ)=F(y^\circ,w^\circ)$.
	\begin{enumerate}
		\item[(a)] If $\hat y$ is obtained from the optimality-MIQP~\eqref{eq:optimality}, then
		\begin{align}
		0&\leq J_{\NLP}(\hat y)-J_{\NLP}(y^\circ)\label{eq:nlp-loss}\\
		&\leq K\sqrt{(N-M)n_v}
		+C\bigl(\norm{\hat y-y^*}^2+\norm{y^\circ-y^*}^2\bigr).\notag
		\end{align}
		\item[(b)] If $\hat y$ is obtained from the feasibility-MIQP~\eqref{eq:feasibility}, then
			\begin{align}
			 J_{\NLP}(\hat y)-J_{\NLP}(y^\circ)\notag
			&\leq C\bigl(\norm{\hat y-y^*}^2+\norm{y^\circ-y^*}^2\bigr)\\
			&-(J_\QP^* - L_M).\label{eq:nlp-feasibility-value}
		\end{align}
	\end{enumerate}
	\end{theorem}
\begin{proof}
	By Assumption~\ref{ass:integer-locality}, Lemma~\ref{lem:value-consistency} applies at $\hat y$ and $y^\circ$.
	Adding and subtracting the corresponding QP values gives
\begin{align}
	&J_{\NLP}(\hat y)-J_{\NLP}(y^\circ)\label{eq:nlp-qp-loss}\\
	&=J_{\NLP}(\hat y)-J_{\QP}(\hat y)
	+J_{\QP}(\hat y)-J_{\QP}(y^\circ)\notag\\
	&\phantom{={}}+J_{\QP}(y^\circ)-J_{\NLP}(y^\circ)\notag\\
	&\leq J_{\QP}(\hat y)-J_{\QP}(y^\circ)
	+|J_{\NLP}(\hat y)-J_{\QP}(\hat y)|\notag\\
	&\phantom{\leq{}}+|J_{\QP}(y^\circ)-J_{\NLP}(y^\circ)|\notag\\
	&\overset{\eqref{eq:value-consistency}}{\leq} J_{\QP}(\hat y)-J_{\QP}(y^\circ)
	+C\bigl(\norm{\hat y-y^*}^2+\norm{y^\circ-y^*}^2\bigr).\notag
\end{align}
	Since $y^\circ\in\mathcal B_\epsilon(y^*)$, QP~\eqref{eq:qp-value} with $y=y^\circ$ admits a solution.
	Let $J_{\QP}^*$ denote the globally optimal value of the full MIQP~\eqref{eq:miqp}.
	Together with Theorem~\ref{th:horizon-gap} and $y^\circ\in\Y$, this yields
	\begin{equation}\label{eq:qp_circ_bound}
		L_M\leq J_{\QP}^*\leq J_{\QP}(y^\circ).
	\end{equation}
	In case~(a), optimality of the reconstruction gives $J_{\QP}(\hat y)=U_M$.
	Comparison with the full MIQP and Theorem~\ref{th:horizon-gap} yields
	\begin{align*}
		J_{\QP}(\hat y)-J_{\QP}(y^\circ)
		&\overset{\eqref{eq:qp_circ_bound}}{\leq} U_M-L_M
		\overset{\eqref{eq:horizon-gap-M}}{\leq} K\sqrt{(N-M)n_v}.
	\end{align*}
	Substitution into~\eqref{eq:nlp-qp-loss} proves the upper bound in~\eqref{eq:nlp-loss}.
	Its lower bound follows from feasibility of the polished pair $(\hat y,\hat w)$ and global optimality of $(y^\circ,w^\circ)$, hence $J_{\NLP}(y^\circ) \leq J_{\NLP}(\hat y)$.

	In case~(b), optimality of the reconstruction gives $J_{\QP}(\hat y)=L_M$.
	Using~\eqref{eq:qp_circ_bound}, we obtain
	\begin{align*}
		J_{\QP}(\hat y)-J_{\QP}(y^\circ) \leq L_M-J_{\QP}^*
		=-(J_{\QP}^*-L_M) \leq 0.
	\end{align*}
	Substitution into~\eqref{eq:nlp-qp-loss} proves~\eqref{eq:nlp-feasibility-value}.
\end{proof}

For the optimality relaxation, the nonlinear loss is bounded by the horizon-tail gap and a local quadratic modeling error.
The feasibility relaxation does not enforce integrality over the full horizon, which tightens the upper bound by the relaxation gap $J_{\QP}^*-L_M$ at the expense of full integer feasibility.
At $M=N$, the upper bounds of \eqref{eq:nlp-loss} and \eqref{eq:nlp-feasibility-value} are identical and reduce to the form of~\cite[Theorem~2]{Buerger2023}.

\section{Advanced-step mixed-integer NMPC}\label{sec:advanced}
In real-time MPC, the delay between receiving a state measurement and applying the control input should be as short as possible.
To reduce this feedback delay, we reorder and reduce the computations of the three-step decomposition  into a preparation phase, which solves the continuous relaxation with a predicted state, and a feedback phase, which solves a horizon-tail relaxed MIQP with the measured state.

\subsection{Preparation and feedback}
Following the advanced-step ideas~\cite{Zavala2009,Nurkanovic2019a,Frey2024a}, we separate computations performed before the next measurement from those needed to compute the feedback.
At sampling time $t_k$, denote the true initial state by $\bar x_0^k$ and its prediction by $\tilde x_0^k$.

During the \textit{preparation phase}, while the previously computed input $u_0^{k-1}$ is applied, we solve the continuous relaxation at $\tilde x_0^k$ to obtain $(y_k^*,w_k^*)$ and its multipliers.
This completes \eqref{step:relaxation} and provides the linearization point for the MIQP.
For the optimality relaxation, we also compute the binary sequence $y_k^{**}$ used to fix the tail.
In principle, this can be obtained with any feasible method, e.g., heuristics, rule based controllers or learned policies~\cite{Boldocky2026}.
In our implementation, we solve the CIA-MILP via \texttt{pycombina} to obtain $y_k^{**}$.
Optionally, the MIQP can be solved at $\tilde x_0^k$ to provide a warm start for the feedback phase MIQP with $\bar x_0^k$.

Once $\bar x_0^k$ is available, the \textit{feedback phase} updates the initial state to $\bar x_0^k$ and solves the MIQP.
All other data in the MIQP remain unchanged.
Algorithm~\ref{alg:advanced-step} summarizes the procedure.
All three MIQP variants enforce integrality of the first control input, including the case $M_o=0$, where the entire binary sequence is prescribed, and only a QP is solved in the feedback phase.

\vspace{-0.3cm}
\begin{algorithm}[h!]
	\caption{Advanced-step mixed-integer nonlinear MPC \\ (AS-MINMPC)}\label{alg:advanced-step}
	\begin{algorithmic}[1]
		\Require Horizons $N,M$ and a choice of~\eqref{eq:miqp}, \eqref{eq:optimality}, or~\eqref{eq:feasibility}
		\For{each sample $k$}
		\Statex \textbf{Preparation before $t_k$}
		\State Predict the next initial state $\tilde x_0^k$.
		\State Solve the relaxed NLP of~\eqref{step:relaxation} with $\tilde x_0^k$ to obtain $(y_k^*,w_k^*)$ and its multipliers.
		\State Set up the MIQP of~\eqref{step:reconstruction} at the lin. point $(y_k^*,w_k^*)$.
		\State For~\eqref{eq:optimality}, compute $y_k^{**}\in\Y$ and fix its tail.
		\State \textcolor{gray}{\emph{Optional:} Solve the selected MIQP at $\tilde x_0^k$ to obtain a warm start.}
		\Statex \textbf{Feedback after measuring $\bar x_0^k$}
		\State Set the initial state in the MIQP to $\bar x_0^k$.
		\State Solve the selected MIQP of~\eqref{step:reconstruction} using the available warm start, to obtain $(\hat y_k,\hat w_k)$.
		\State \textcolor{gray}{\emph{Optional:} Solve the polishing NLP in~\eqref{step:polishing} with $\hat y_k$.}
		\State Apply the first continuous and integer control inputs.
		\State Shift the solution for the next preparation phase.
		\EndFor
	\end{algorithmic}
	\end{algorithm}
	\vspace{-0.5cm}
\subsection{MIQP step as a full SQP step}\label{sec:sqp-step}
To reduce feedback delay, the proposed algorithm does not in general solve the polishing NLP in~\eqref{step:polishing} to convergence.
We observe that the feedback MIQP simultaneously reconstructs the integer decisions $\hat y_k$ and performs one SQP step toward the corresponding polishing solution, yielding $\hat w_k\approx w_{\NLP}(\hat y_k)$.

For the full and optimality-relaxed MIQPs, fixing the reconstructed integers yields the SQP subproblem obtained by linearizing the NLP in $(y,w)$ at $(y^*,w^*)$ and imposing $y=\hat y$.
For the feasibility relaxation, the same interpretation holds after fixing only the binary prefix and retaining the tail as continuous variables.

This interpretation connects the method to real-time iteration (RTI) and advanced-step RTI~\cite{Diehl2009c,Nurkanovic2019a,Frey2024a} for continuous NLPs.
If the integer decisions are the same for the measured and predicted states, i.e., $y(\bar x_0^k)=y(\tilde x_0^k)$, we expect a local predictor--corrector error analysis for the proposed AS-MINMPC method to yield bounds similar to~\cite{Frey2024a}.
If $y(\bar x_0^k)\neq y(\tilde x_0^k)$, the solution map may be discontinuous, similarly to nonsmooth MPC~\cite{Nurkanovic2026b}, which would need to be accounted for in the error analysis.
Due to page restrictions, we leave a detailed error analysis for future work.

\section{Numerical example}\label{sec:example}
We compare several variants of the proposed method with the state-of-the-art MINLP solver \texttt{CAMINO}~\cite{Ghezzi2026a}, and a rule-based controller (RBC) from~\cite{Boldocky2026} on a parallel-chiller example~\cite{Boldocky2026}.

\subsection{Implementation}
We implement several variants of Algorithm~\ref{alg:advanced-step} (including the optional steps, and solving MILPs instead of MIQPs) in the accompanying repository.\footnote{\url{https://github.com/nosnoc/advanced-minmpc}}
The implementation uses \texttt{IPOPT}~\cite{Waechter2006} for NLPs, \texttt{BnB-DAQP}~\cite{Arnstrom2023} for MIQPs via \texttt{CasADi}~3.8.0~\cite{Andersson2019}, and \texttt{pycombina}~\cite{Buerger2020a} for integer-feasible reconstruction.
Thus, the proposed method uses only open-source software packages.

\subsection{Benchmark example}
To evaluate the proposed algorithms, we consider a data center parallel-chiller plant example from~\cite[Sec.~2]{Boldocky2026}.
The objective is to meet the cooling demand while reducing power consumption and chiller switching frequency.
For $n_v$ chillers, the state vector is $x=(T_r,T_s^{(1)},\ldots,T_s^{(n_v)})$, comprising the shared return and individual supply temperatures.
The continuous controls are $u=(T_e^{(1)},\ldots,T_e^{(n_v)},\dot m^{(1)},\ldots,\dot m^{(n_v)})$, comprising the evaporator temperatures and mass flows, while $v=(v^{(1)},\ldots,v^{(n_v)})$ contains the binary on/off decisions for each chiller operating in parallel.
The cooling powers and thermal dynamics are:
	\begin{align*}
		Q^{(i)}(t)&=\eta_r c_p\dot m^{(i)}v^{(i)}(T_r(t)-T_s^{(i)}(t)),\\
		\dot T_r(t)&=\frac{1}{C_r}(\widetilde Q_{\mathrm{load}}(t)-\sum_{i=1}^{n_v}Q^{(i)}(t)),\\
		\dot T_s^{(i)}(t)&=-\frac{c_p\eta_s}{C^{(i)}}v^{(i)}(t)\dot m^{(i)}(t)(T_s^{(i)}(t)-T_e^{(i)}(t)).
	\end{align*}
The chiller and pump power consumption are:
	\begin{align*}
		P_{\mathrm{chiller}}^{(i)}(t)
		&=\frac{Q^{(i)} (t)}{\mathrm{COP}^{(i)} (t)}+\rho^{(i)}v^{(i)}(t),\\
		\mathrm{COP}^{(i)}(t)\!&=\!a_0^{(i)}+a_1^{(i)}\frac{Q^{(i)}(t)}{Q_{\max}^{(i)}}
		+a_2^{(i)}\left(\frac{Q^{(i)}(t)}{Q_{\max}^{(i)}}\right)^2\!,\\
		P_{\mathrm{pump}}^{(i)}(t)	& =\gamma^{(i)} v^{(i)}(t)\left(\dot m^{(i)}(t) \right)^3.
	\end{align*}
All model parameters, bounds, and load profiles are taken from~\cite[Table 1]{Boldocky2026}.
Using the compact notation of the OCP~\eqref{eq:ocp}, we have
\begin{mini!}|s|[2]
	{x,u,v}{
		\sum_{k=0}^{N-1}\ell_k(x_k,u_k,v_k)
		+\sum_{k=0}^{N-1}\norm{v_k-v_{k-1}}_R^2
		\label{eq:chiller-ocp-objective}}
	{\label{eq:chiller-ocp}}{}
	\addConstraint{x_0}{=\bar x_0,\label{eq:chiller-ocp-initial}}
	\addConstraint{x_{k+1}}{=f_k(x_k,u_k,v_k),\label{eq:chiller-ocp-dynamics}}{\quad k\!=\!0,\ldots,N-1}
	\addConstraint{x_{\min}\leq x_k}{\leq x_{\max},\label{eq:chiller-ocp-state}}{\quad k=0,\ldots,N}
	\addConstraint{u_{\min}\leq u_k}{\leq u_{\max},\label{eq:chiller-ocp-control}}{\quad k\!=\!0,\ldots,N-1}
	\addConstraint{0\leq Q_k^{(i)}}{\leq Q_{\max}^{(i)},\label{eq:chiller-ocp-capacity}}{\quad {k=0,\ldots,N-1}}
	\addConstraint{\textstyle\sum_{i=1}^{n_v}v_k^{(i)}}{\geq1\label{eq:chiller-ocp-active},}{\quad k=0,\ldots,N-1}
	\end{mini!}
Here, $v_{-1}$ denotes the previously applied binary control input.
The stage costs terms $\ell_k$ consist of the power consumption and load demand tracking error:
$\ell_k(x_k,u_k,v_k)
	:=\sum_{i=1}^{n_v}\left(P_{\mathrm{chiller},k}^{(i)}
	+P_{\mathrm{pump},k}^{(i)}\right)+\lambda_W\left(\sum_{i=1}^{n_v}Q_k^{(i)}
	-Q_{\mathrm{demand},k}\right)^2$.
%
We use the cost weights $R=200I$ and $\lambda_W=0.01$.
The discrete-time map $f_k$ incorporates the prescribed effective demand and one fourth-order Runge--Kutta (RK4) step per control interval of $180\,\mathrm{s}$.
Closed-loop simulation uses four RK4 substeps per control interval to obtain the next measured state.

\begin{figure}[t]
	\centering
	\includegraphics[width=1.00\columnwidth]{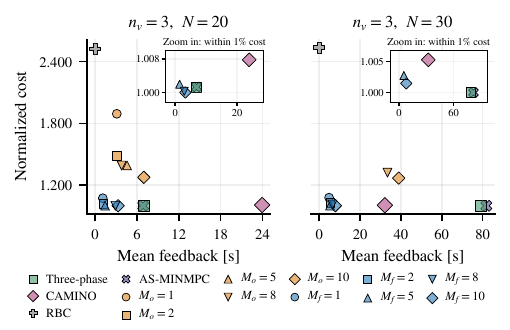}
	\vspace{-0.8cm}
	\caption{Mean feedback time and normalized realized closed-loop cost for three chillers with $N=20$ and $N=30$.
	Horizon-tail relaxations yield greater AS-MINMPC speedups for longer prediction horizons.}
	\label{fig:chiller-pareto}
\end{figure}

\begin{figure}[t]
	\centering
	\vspace{-0.4cm}
	\includegraphics[width=0.95\columnwidth]{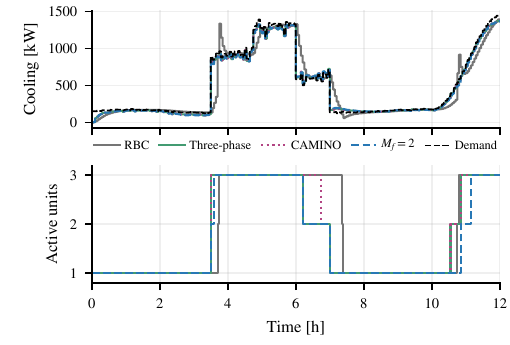}
	\vspace{-0.45cm}
	\caption{Representative MPC trajectories for delivered cooling and demand (top), and number of active chillers (bottom), for $N =30$.}
	\label{fig:chiller-trajectories}
	\vspace{-0.7cm}
\end{figure}

All controllers use the same load profile (Fig.~\ref{fig:chiller-trajectories}) and initial states over $240$ MPC steps (12 h of sim. time).
Preparation times are never the computational bottleneck, therefore, we focus on and report feedback times, which are dominated by the MIQP time.
The realized closed-loop cost sums the switching and stage costs $\ell_k$ incurred by the applied inputs over the simulation and is normalized by the lowest cost among all tested algorithms for each experiment.
We compare several variants of the proposed AS-MINMPC method with different $M_f$ and $M_o$ values, the three-phase algorithm~\cite{Buerger2023} of Sec.~\ref{sec:decomposition}, and a state-of-the-art MINLP method: \texttt{SB-MIQP}, implemented in \texttt{CAMINO}~\cite{Ghezzi2026a}.
The latter two perform all computations in the feedback phase.
We use the default settings of \texttt{SB-MIQP}, except for a time limit of $30\,\mathrm{s}$ per MPC step, with \texttt{Gurobi} as MIQP solver.
A single full SB-MIQP iteration is equivalent to the three phase algorithm~\cite{Ghezzi2026a}.

Fig.~\ref{fig:chiller-pareto} compares mean feedback times and normalized costs for $N=20,30$ with $n_v=3$ chillers.
For $N=20$, the variant $M_f=5$ computes feedback $16.68$ times faster than \texttt{CAMINO}, with a slightly lower closed-loop cost.
For $N=30$, it reduces the mean and maximum observed feedback times by factors of $15.56$ and $28.30$, respectively, relative to full AS-MINMPC, with only a $0.27\%$ increase in closed-loop cost.
RBC is computationally cheaper, but its cost is approximately $2.5$ times that of the best optimization-based controller.
Figure~\ref{fig:chiller-trajectories} shows that $M_f=2$ closely follows the cooling and chiller activation trajectories of the three-phase method, with better demand tracking than RBC.

In this example, feasibility-relaxed MIQPs yielded substantial speedups, often with minor losses of optimality below 1\%.
Optimality-relaxed MIQPs resulted in higher costs in these experiments due to the less favorable integer-sequences obtained with \texttt{pycombina}.
Our experiments on other examples available in the public repository showed that this relative performance is problem dependent.
But in general, horizon-tail approximations substantially reduced feedback delay while preserving near-reference control performance across all implemented examples.

\color{black}

\section{Conclusions}
We presented a real-time mixed-integer nonlinear MPC method that reorders the computations of a decomposition algorithm~\cite{Sager2009,Buerger2023} into preparation and feedback phases.
The feedback MIQP reconstructs the integer controls and simultaneously performs one SQP step toward the solution of the problem with the latest state estimate.
Horizon-tail approximations in the MIQP reduce the combinatorial complexity associated with the horizon length.
Numerical results show reductions of more than $95\%$ in the maximum observed feedback time, with only minor changes in closed-loop cost.
The computational advantage of the proposed real-time algorithms becomes increasingly pronounced as the prediction horizon grows.
The feedback phase still requires solving an MIQP, which can be expensive when many integer control variables are present.
This motivates development of more efficient open-source MIQP solvers.
\bibliographystyle{IEEEtran}

\end{document}